\documentclass[12pt]{article}
\usepackage{amssymb,amsfonts,amsmath,latexsym,epsf,url}
\usepackage[usenames,dvipsnames]{pstricks}
\usepackage{pstricks-add}
\usepackage{subfigure}
\usepackage{graphicx}
\usepackage{color,soul}
\usepackage{epsfig}
\usepackage{tikz}
\usepackage[colorlinks]{hyperref}
\usepackage{cite}
\usepackage{algorithm}
\usepackage{algpseudocode}
\usepackage{graphicx}
\usepackage{booktabs}
\usepackage{multirow}

\newtheorem{theorem}{Theorem}[section]

\newtheorem{proposition}[theorem]{Proposition}
\newtheorem{observation}[theorem]{Observation}

\newtheorem{lemma}[theorem]{Lemma}

\newcommand{\qed}{\hfill $\square$\medskip}

\begin{document}

\def\nt{\noindent}

\title{On the spectrum of two families of non-distance-regular graphs}

\author{
	Ali Zafari$^1$\footnote{Corresponding author} \and	
Saeid Alikhani$^{2}$
}

%\date{\today}

\maketitle

\begin{center}

$^1$Department of Mathematics, Faculty of Science,
Payame Noor University, P.O. Box 19395-4697, Tehran, Iran\\ 
{\tt zafari.math@pnu.ac.ir}
\medskip

$^{2}$Department of Mathematical Sciences, Yazd University, 89195-741, Yazd, Iran\\
{\tt alikhani@yazd.ac.ir}

\end{center}
	
%%%%%%%%%%%%%%%%%%%%%%%%%%%%%%%%%%%%%%%%%%%%%%%%%%%%%%%%%%%%%%%%%%%%%%%%%%%%%%%%%
%%%%%%%%%%%%%%%%%%%%%%%%%%%%%%%%%%%%%%%%%%%%%%%%%%%%%%%%%%%%%%%%%%%%%%%%%%%%%%%%%
\begin{abstract}
This paper addresses the challenge of spectral analysis and structural investigation for graphs that are not distance-regular, where computing the spectrum using standard methods based on equitable and orbit partitions can be complex. 
Our main objective is to determine all eigenvalues of the extended graph $E(2.O_k)$ by leveraging the relationship between its equitable and orbit partitions. While the integral nature of this graph has been previously studied, we introduce a novel approach to demonstrate the utility of this method in finding the complete set of distinct eigenvalues for a class of non-distance-regular graphs.
Specifically, we first establish that $E(2.O_k)$ is a vertex-transitive graph with diameter $k$, contrasting with the diameter of $2.O_k$, which is $2k-1$. We also determine the automorphism group of $E(2.O_k)$ and prove that it is an integral graph, meaning all eigenvalues of its adjacency matrix are integers. A significant result is the determination of the multiplicity for all distinct eigenvalues of $E(2.O_k)$.
Additionally, we extend our method to the enhanced Johnson graph $EJ(2m,m)$. Although its eigenvalues are known from prior work, the multiplicity of these distinct eigenvalues has not yet been calculated. We use our techniques to fully determine the multiplicity of all distinct eigenvalues for $EJ(2m,m)$.
\end{abstract}


\medskip
\noindent{\bf AMS Subj.\ Class.:}  05C25, 94C15.

%%%%%%%%%%%%%%%%%%%%%%%%%%%%%%%%%%%%%%%%%%%%%%%%%%%%%%%%%%%%%%%%%%%%%%%%%%%%%%%%%
%%%%%%%%%%%%%%%%%%%%%%%%%%%%%%%%%%%%%%%%%%%%%%%%%%%%%%%%%%%%%%%%%%%%%%%%%%%%%%%%%

\section{Introduction}
\label{sec:introduction}

The {Odd graph $O_k$}, first introduced by Balaban, Farcussiu, and Banica \cite{paper1.2} and Biggs \cite{paper3.1}, is a well-studied family of graphs \cite{paper5,paper4,paper7}. For $k \geq 3$, the vertices of $O_k$ are the $(k-1)$-subsets of a set $\Omega$ of size $2k-1$. Two vertices (subsets) are adjacent if they are disjoint.

The related graph, the {bipartite double $2.O_k$}, has $2\binom{2k-1}{k-1}$ vertices, represented as pairs $(v, i)$ where $v \subset \Omega$ is a $(k-1)$-subset and $i \in \{0, 1\}$. Adjacency requires that the subsets are disjoint ($v \cap w = \emptyset$) and the indices differ ($i \neq j$). It is known that $2.O_k$ is distance transitive (and thus distance regular) with valency $k$ and diameter $2k-1$ \cite{paper3}. Crucially, $2.O_k$ is isomorphic to the middle cube graph $MQ_k$ \cite{paper6a}.

The primary object of study is the {extended graph $E(2.O_k)$}, which shares the same vertex set as $2.O_k$ but includes additional cross edges $E_2$ connecting a vertex $(v, i)$ to its complement-indexed vertex $(v, i^c)$, where $i^c$ is the complement index ($1^c=0, 0^c=1$). This addition makes $E(2.O_k)$ a regular bipartite graph of degree $k+1$.

Let $n,m$ be two positive integers with, $n>m>1$.  The Johnson graph $J(n,m)$ is defined as the graph whose vertex set is $V=\{v\mid v\subseteq I=\{1,\ldots,n\},|v|=m\}$, where two vertices $u$, $v$ are adjacent if and only if $|u\cap v|=m-1$. It is well known that $J(n,m)$ is a vertex transitive graph \cite{paper7}, and  belongs to a symmetric association scheme known as the Johnson scheme, see \cite{bbb2,paper6}.
For the Johnson scheme $J(n,m)$, the Bose–Mesner algebra $\mathcal{A}$ is the complex algebra generated by distance adjacency matrices $A_0, A_1, \dots, A_m$.
The paper also investigates the {enhanced Johnson graph $EJ(2m,m)$}, defined in \cite{w-3}. It is constructed from the standard Johnson graph $J(2m,m)$ by adding edges between complementary subsets ($v$ and $v^c$). While the integrality of $EJ(2m,m)$ is known \cite{w-3}, the multiplicities of its eigenvalues have not yet been determined.

The main goal of this work is to determine the complete spectrum of $E(2.O_k)$, including all distinct eigenvalues and their multiplicities, by using the relationship between equitable partitions and orbit partitions. Furthermore, we aim to:
\begin{itemize}
	\item Establish the structural properties of $E(2.O_k)$, demonstrating that it is {vertex transitive} but {not distance regular}.
	\item Explicitly determine the {automorphism group} of $E(2.O_k)$, showing it is isomorphic to $\mathbb{Z}_2 \times \text{Sym}(2k-1)$.
	\item Compute the multiplicities of all distinct eigenvalues for the enhanced Johnson graph $EJ(2m,m)$.
\end{itemize}

%%%%%%%%%%%%%%%%%%%%%%%%%%%%%%%%%%%%%%%%%%%%%%%%%%%%%%%%%%%%%%%%%%%%%%%%%%%%%%%%%
%%%%%%%%%%%%%%%%%%%%%%%%%%%%%%%%%%%%%%%%%%%%%%%%%%%%%%%%%%%%%%%%%%%%%%%%%%%%%%%%%

\section{Preliminaries} 
This section introduces the fundamental concepts, notation, and prior results from graph theory and algebraic graph theory that are essential for understanding the main results presented in this paper.

	The {diameter} of a connected graph is the maximum distance between any two vertices. The {girth} is the length of the shortest cycle in the graph.
	An {automorphism} of a graph $\Gamma$ is a permutation of its vertices that preserves adjacency. The set of all automorphisms forms the {automorphism group}, denoted $\text{Aut}(\Gamma)$.
	A graph $\Gamma$ is {vertex transitive} if $\text{Aut}(\Gamma)$ acts transitively on the vertex set $V(\Gamma)$ (i.e., there is only one orbit). The graph $\Gamma$ is {distance transitive} if for any two pairs of vertices $(u, v)$ and $(x, y)$ such that the distance $d(u, v) = d(x, y)$, there exists an automorphism mapping $u$ to $x$ and $v$ to $y$.

	A partition $\pi = \{C_1, \ldots, C_r\}$ of the vertex set $V(\Gamma)$ is {equitable} if, for any two cells $C_i$ and $C_j$, the number of neighbors that a vertex in $C_i$ has in $C_j$ is a constant value $b_{ij}$, independent of the specific vertex chosen in $C_i$.

	Let $H \leq \text{Aut}(\Gamma)$ be a subgroup. The partition of $V(\Gamma)$ defined by the orbits of $H$ is called an {orbit partition}. The {quotient graph $\Gamma/\pi$} is a directed graph whose vertices are the cells of an equitable partition $\pi$, with adjacency determined by the $b_{ij}$ constants. The matrix $B_\pi = (b_{ij})$ is the adjacency matrix of the quotient graph.

\begin{proposition}{\rm\cite{paper7,w-2}}
	\label{prop:orbit_equitable}
	An orbit partition is an equitable partition.
\end{proposition}

	A regular graph $\Gamma$ with diameter $d$ is {distance regular} if, for any two vertices $u$ and $v$ at distance $r$, the following intersection numbers are constant:
	\begin{itemize}
		\item $c_r = |\Gamma_{r-1}(v) \cap \Gamma_1(u)|$
		\item $b_r = |\Gamma_{r+1}(v) \cap \Gamma_1(u)|$
	\end{itemize}
	The partition of vertices based on distance from a fixed vertex $v$, denoted $\Gamma_r(v) = \{\Gamma_0(v), \ldots, \Gamma_d(v)\}$, is called the {distance partition}.

\begin{theorem}{\rm\cite{paper1E-1}}\label{a.3}
	\label{thm:dist_reg_equitable}
	For any graph $\Gamma$, the distance partition from any vertex $v$ is equitable if and only if $\Gamma$ is a distance regular graph.
\end{theorem}

\begin{observation}{\rm\cite{paper7}}
	\label{obs:dist_trans_orbit}
	If $\Gamma$ is distance transitive, then for every vertex $v$, the distance partition $\{\Gamma_0(v), \ldots, \Gamma_d(v)\}$ is an orbit partition under the action of the vertex stabilizer group $\text({Aut}(\Gamma))_v$.
\end{observation}

	A homomorphism $f: V(\Gamma) \to V(\Lambda)$ is a {covering map} if it is a surjective local isomorphism (i.e., for every vertex $w$ in $\Lambda$, the mapping from the neighbors of a vertex in the fiber $f^{-1}(w)$ to the neighbors of $w$ is bijective). If $f$ exists, $\Gamma$ is said to {cover} $\Lambda$.

\begin{theorem}{\rm\cite{paper7}}
	\label{thm:quotient_eigenvalues}
	Let $\Gamma$ be a graph with an equitable partition $\pi$. Every eigenvalue of the quotient matrix $B_\pi$ is also an eigenvalue of the adjacency matrix $A$ of $\Gamma$.
\end{theorem}

\begin{theorem}{\rm\cite{paper7}}
	\label{thm:singleton_eigenvalues}\label{b.6}
	Let $\Gamma$ be a vertex transitive graph, and $\pi$ an orbit partition of a subgroup $H \leq \text{Aut}(\Gamma)$. If $\pi$ includes a singleton cell $\{u\}$, then every eigenvalue of $\Gamma$ is an eigenvalue of the quotient graph $\Gamma/\pi$.
\end{theorem}

%%%%%%%%%%%%%%%%%%%%%%%%%%%%%%%%%%%%%%%%%%%%%%%%%%%%%%%%%%%%%%%%%%%%%%%%%%%%%%%%%
%%%%%%%%%%%%%%%%%%%%%%%%%%%%%%%%%%%%%%%%%%%%%%%%%%%%%%%%%%%%%%%%%%%%%%%%%%%%%%%%%

\section{Main Results}
In this section we study some algebraic aspects of $E(2.O_k)$ and also investigate spectrum of $E(2.O_k)$ and enhanced Johnson graph $EJ(2m,m)$.

\subsection{Some algebraic aspects of $E(2.O_k)$ }

\begin{proposition} \label{c.1}
	The graph $E(2.O_k)$ is a vertex transitive graph.
\end{proposition} 
	\begin{proof}
		Let $\Lambda=E(2.O_k)$ and  $\Omega = \{1, 2, ... , 2k-1\}$. It is easy to prove that
		the graph $\Lambda$ is a regular bipartite graph of degree $k + 1$. Now, if
		\[
		V_1=\{ (v, i)  \,\ | \,\, \text{$v$ is a $(k-1)$-subset of $\Omega$}; i\in\{ 0, 1\} \},
		\] 
		and
		\[
		V_2=\{ (v, i^c	)  \,\ | \,\, \text{$v$ is a $(k-1)$-subset of $\Omega$}; i^c\in\{ 0, 1\} \},
		\]
		then $ V(\Lambda)=V_1\cup V_2$, such that every edge of $\Lambda$ has one end in $V_1$ and one end in $V_2$, and hence
		$| V_1 |=| V_2 |=\binom{2k-1}{k-1}$.
		By the following steps for each pair of distinct vertices from  $V(\Lambda)$, we show that $\Lambda$ is a vertex transitive graph.
		
		(i)
		If both vertices $(u,i)$ and $(v, i)$ lie in $V_1$ and $|u \cap v| = t$, where $0\leq t\leq k-2$, then we may assume
		$u =\{ x_1, ... , x_t, u_1, ... , u_{(k-1)-t}\}$ and
		$v =\{ x_1, ... , x_t, v_1, ... , v_{(k-1)-t}\}$, where $x_i, u_i, v_i \in \Omega $. Let $\sigma$ be a permutation of $Sym(\Omega )$, such that $\sigma(x_i)=x_i$, $\sigma(u_i)=v_i$ and $\sigma(w_i)=w_i$, where $w_i\in \Omega -(u\cup v)$.
		So, $\sigma$ induces an automorphism $\widetilde{\sigma}:V(\Lambda)\rightarrow V(\Lambda)$ by
		$$\widetilde{\sigma}(\{x_1, ... , x_t, u_1, ... , u_{(k-1)-t}\}, i)=(\{ \sigma (x_1), ... , \sigma (x_t), \sigma( u_1), ... , \sigma (u_{(k-1)-t})\}, i).$$ Therefore, $\widetilde{\sigma}(u, i)=(v, i).$
		
		(ii)
		We define the mapping $\theta: V(\Lambda)\rightarrow V(\Lambda) $ by  $\theta(v, i) = (v, i^c)$ for every $(v, i)$ in $V(\Lambda)$. It is easy to prove that $\theta$ is an automorphism of $\Lambda$. So, if both  vertices $(u, i^c)$ and $(v, i^c)$ lie in $V_2$, then $\theta(u, i^c)\in V_1$, and $\theta(v, i^c)\in V_1$. Therefore, there is
		an automorphism $\widetilde{\sigma}$ in $Aut(\Lambda)$ such that, $\widetilde{\sigma}(\theta(u, i^c))=\theta(v, i^c)$. Thus, 
		$(\theta^{-1}\widetilde{\sigma}\theta)(u, i^c)=(v, i^c)$.
		
		(iii)
		Now, let $(u, i)\in V_1$ and $(v, i^c)\in V_2$, so $\theta(v, i^c)\in V_1$, and so, there is an automorphism $\widetilde{\sigma}$
		in $Aut(\Lambda)$ such that $\widetilde{\sigma}(u, i)=\theta(v, i^c)$. Thus,
		$(\theta^{-1}\widetilde{\sigma})(u, i)=(v, i^c)$.\qed
		%\end{itemize}
	\end{proof}

%%%%%%%%%%%%%%%%%%%%%%%%%%%%%%%%%%%%%%%%%%%%%%%%%%%%%%%%%%%%%%%%%%%%%%%%%%%%%%%%%
\begin{proposition} \label{c.2}
	The diameter of $E(2.O_k)$ is $k$, and its girth is $4$.
\end{proposition}
	\begin{proof}
		Let $\Lambda=E(2.O_k)$ and let $(u,i)$ and $(w,j)$ be any two vertices in $\Lambda$. Based on \cite{paper1Ee-1}, we can verify that the diameter of $E(2.O_k)$ is $k$. Now, by another approach we show that the diameter of $E(2.O_k)$ is $k$.
		If the distance between two distinct vertices $u$ and $w$ in Odd graph $O_k$ is equal to $d$, that is
		$d_{O_k}(u,w)=d$ then $d\leq k-1$, because it is well known that the diameter $O_k$ is equal to $k-1$. If we consider a path $P_{O_k}:u=u_0, u_1, ..., u_d=w$  in  Odd graph $O_k$, then by lift this path to $\Lambda$ and using cross edges
		$P_{\Lambda}:(u,a_0)=(u_0,a_0), (u_1,a_1), (u_2,a_2) ..., (u_d,a_d)=(w,a_d)$, where $a_0=i$ and each step requires $a_{t+1}=a^{c}_t$, we have 
		$a_d=i$ if $d$ is even, and  $a_d=i^c$ if $d$ is odd. Thus, if $d=k-1$ and 
		$a_d\neq j$, then the distance between $(u,i)$ and $(w,j)$ is exactly $k$. Also, the girth of $E(2.O_k)$ is equal to $4$, because if two distinct vertices
		$(u, i)$ and $(w, i^c)$, are adjacent in $E(2.O_k)$, where $u \cap w =\emptyset$, then $(u, i^c)$ and $(w, i)$, are also adjacent in $E(2.O_k)$. Hence, any edge of 
		$E(2.O_k)$ lie in a cycle of length $4$, because $E(2.O_k)$ is a vertex transitive graph. \qed
	\end{proof}
%%%%%%%%%%%%%%%%%%%%%%%%%%%%%%%%%%%%%%%%%%%%%%%%%%%%%%%%%%%%%%%%%%%%%%%%%%%%%%%%%

\begin{proposition} \label{c.3}
	The graph $E(2.O_k)$ cannot be a  distance regular graph.
\end{proposition} 
	\begin{proof}
		For each vertex $(v, i)$ of $E(2.O_k)$, there is a unique vertex $(v, i^c)$ in $E(2.O_k)$ so that two vertices $(v, i)$ and $(v, i^c)$ are adjacent in $E(2.O_k)$, because the bipartite double graph  $2.O_k$ of the Odd graph $O_k$ has the property that, for each vertex $(v, i)$, there is a unique vertex in $2.O_k$  at distance $2k-1$ from $(v, i)$, say $(v, i^c)$. Now, let $\Lambda=E(2.O_k)$ and let $\Lambda_r((v, i))=\{\Lambda_0(v, i), ... , \Lambda_{k}(v, i)\}$ be the distance partition of $E(2.O_k)$ for each vertex $(v, i)$ of $E(2.O_k)$, where $r$ is a non-negative integer not exceeding $k$, the diameter of $E(2.O_k)$. Hence the vertex $(v, i^c)$ lie in the cell $\Lambda_1(v, i)$ of  the distance partition $\Lambda_r((v, i))$. 
		Also, for two distinct vertices $(w, i^c)$  and $(v, i^c)$ in the cell $\Lambda_1(v, i)$ there are distinct vertices  $(x, i)$ and $(w, i)$ in the cell $\Lambda_2(v, i)$ so that $(w, i^c)$ is adjacent to $(x, i)$, $(w, i)$ and $(v, i^c)$ is adjacent to  $(w, i)$. Therefore,
		$d((v,i), (x, i))=d((v,i), (w, i))=2$, and hence we have $1=|\Lambda_{1}(v, i)\cap \Lambda_{1}(x, i)|\neq|\Lambda_{1}(v, i)\cap \Lambda_{1}(w, i)|=2$.
		Thus, $E(2.O_k)$ cannot be a  distance regular graph. \qed
	\end{proof}
%%%%%%%%%%%%%%%%%%%%%%%%%%%%%%%%%%%%%%%%%%%%%%%%%%%%%%%%%%%%%%%%%%%%%%%%%%%%%%%%%

\begin{theorem} \label{c.4}
	The automorphism group of the graph $\Lambda=E(2.O_k)$
	is identical to the automorphism group of  $\Gamma=2.O_k$.
\end{theorem} 
	\begin{proof}
		Let $H=Aut(\Gamma)$, and $\Omega=\{1, ... , 2k-1\}$. We know that	
		$H\cong \mathbb{Z}_2\times Sym(\Omega)$, (see  ~\cite{paper6} Page 260). Moreover, let	 $G=Aut(\Lambda)$.	
		It is not hard to see that	 $G\leq H$, because if  $f\in G$ and $e=\{(v,i),(u,j)\}$, where $v \cap u =\emptyset$, and $i\neq j$  is an edge in $\Gamma$, then 
		$e$ is an edge in $\Lambda$. In particular any edge of $\Lambda$ as to form $e=\{(v,i),(v,j)\}$, $i\neq j$ or $\{(w, i),(z, j)\}$  so that  $w \cap z =\emptyset$ and $i\neq j$.
		On the other hand, any edge of $\Lambda$  as to form $e=\{(v,i),(v,j)\}$,  $i\neq j$ lie in many $4$ cycles in  $\Lambda$. Also, none of two distinct edges of $\Lambda$ as to form
		$\{(v, i),(u, j)\}$ and $\{(w, i),(z, j)\}$  so that  $v \cap u =\emptyset$, $w \cap z =\emptyset$ and $i \neq j$ are not contained in any cycle of length $4$, where $v$, $u$, $w$ and $z$ are distinct. In particular, two distinct edges of $\Lambda$ as to form
		$\{(v, i),(u, j)\}$ and $\{(v, i),(z, j)\}$  so that  $v \cap u =\emptyset$, $v \cap z =\emptyset$ and $i \neq j$ are not contained in any cycle of length $4$.
		Especially, we  know that the automorphisms preserve cycles, neighbours, hence any  $f\in G$ must map the edge $\{(v, i),(u, j)\}$, $v \cap u =\emptyset$,  $i \neq j$ to edge 
		$\{(w, i),(z, j)\}$, $w \cap z =\emptyset$ and $i \neq j$, and any $f\in G$ must map the edge $\{(v, i),(u, j)\}$, $v \cap u =\emptyset$,  $i \neq j$ to edge 
		$\{(v, i),(z, j)\}$  so that  $v \cap z =\emptyset$ and $i \neq j$.
		Thus,
		for any two edges of $\Lambda$  as  to  form $\{(v, i),(v, j)\}$, $i\neq j$  and $\{(w, i),(z, j)\}$, $w \cap z =\emptyset$,  $i \neq j$
		there is no automorphism from $G$ that takes $\{(v, i),(v, j)\}$ to $\{(w, i),(z, j)\}$. Hence $f(e)=\{f(v,i),f(u,j)\}$ is an edge in $\Gamma$, and hence   $f\in H$. So $G\leq H$. 	
		Now, if we find a subgroup $M$ of $G$ of order $2(2k-1)!$, then we can conclude
		that $G = M=H$.
		Let $K$ be the group
		that is generated by $\theta$, where $\theta: V (\Lambda)\longrightarrow V (\Lambda)$, by  $\theta(v, i)=(v, i^c)$ for every $(v, i)\in V(\Lambda)$, is an automorphism of order $2$ in the graph $\Lambda$.
		Moreover,
		let $L=\{f_{\sigma} \,|\,\,\sigma\in Sym(\Omega) \}$, where $f_{\sigma}: V (\Lambda)\longrightarrow V (\Lambda)$, by $f_{\sigma}(\{x_1, x_2, ... , x_{k-1}\}, i)=(\{\sigma(x_1), \sigma(x_2), ... , \sigma(x_{k-1})\}, i)$ be an automorphism of  $\Lambda$ for every $v\in V(\Lambda)$. By similar way is done in the proof of Theorem 3.6, in \cite{w-1}, and a few changes we can verify that $K$ and $L$ are normal subgroups of $G$ and $K\cap L=1$. Thus, $M=K\times L=G\cong \mathbb{Z}_2 \times Sym(\Omega)=\mathbb{Z}_2 \times Sym(2k-1)$. \qed
	\end{proof}
%%%%%%%%%%%%%%%%%%%%%%%%%%%%%%%%%%%%%%%%%%%%%%%%%%%%%%%%%%%%%%%%%%%%%%%%%%%%%%%%%
%%%%%%%%%%%%%%%%%%%%%%%%%%%%%%%%%%%%%%%%%%%%%%%%%%%%%%%%%%%%%%%%%%%%%%%%%%%%%%%%%
\subsection{Spectrum of  $E(2.O_k)$ and $EJ(2m,m)$}

The adjacency matrix $A=A(\Gamma)$ of a (simple, undirected) graph $\Gamma$ is the $n\times n$ symmetric matrix whose rows and columns are indexed by the vertices of $\Gamma$, and where $A_{xy} = 1$ if and only if  $x$ is adjacent to $y$ (that is $xy\in E$), and $A_{xy} = 0$ otherwise.
If all the eigenvalues of the adjacency matrix of a graph $\Gamma$ are integers,
then we say that $\Gamma$ is an integral graph. The notion of integral graphs was
first introduced by F. Harary and A. J. Schwenk in 1974  ~\cite{paper8}.  In general, the problem of
characterizing integral graphs seems to be very difficult, especially if a graph  cannot be a distance regular graph.
There are some of the recent works in this scope of research
in algebraic graph theory include  \cite{paper0,paper1,paper2}, and \cite{w-2,w-3}.
The square matrix $A(\Gamma)$ is called diagonalizable if it is similar to a diagonal matrix, that is, if there
exists an invertible matrix $P$ such that $P ^{-1}A(\Gamma)P$ is a diagonal matrix.
The eigenvalues of a graph $\Gamma$ are the eigenvalues of the adjacency matrix of $\Gamma$.
The characteristic polynomial of $\Gamma$ with respect to the adjacency matrix $A(\Gamma)$ is the polynomial
$P(\Gamma) = P(\Gamma, \lambda) = det(\lambda I_n - A(\Gamma))$, where $I_n$
denotes the $n \times n$ identity matrix. Because $A(\Gamma)$ is real and symmetric, its eigenvalues are real numbers.
%%%%%%%%%%%%%%%%%%%%%%%%%%%%%%%%%%%%%%%%%%%%%%%%%%%%%%%%%%%%%%%%%%%%%%%%%%%%%%%%%
The spectrum of a graph $\Gamma$  is the list of the eigenvalues of the adjacency matrix of $\Gamma$ together with their multiplicities, and it is denoted by $Spec(\Gamma)$. 
Also,  we denoted  the multiplicity of each  eigenvalue $\lambda$  of the adjacency matrix of $\Gamma$ by $m(\lambda)$. \\
In this section, by using the relationship between equitable partition and orbit partition of graphs we show that $E(2.O_k)$ is an integral graph, although $E(2.O_k)$ cannot be a distance regular graph. Especially, we determine the multiplicity of all distinct eigenvalues of the extended graph $E(2.O_k)$. Furthermore, if we consider the enhanced Johnson graph $EJ(2m,m)$ is defined already, its eigenvalues have been computed in \cite{w-3}; however, the multiplicity of the eigenvalues of the enhanced Johnson graph $EJ(2m,m)$ has not been calculated. Thus, in this paper, we determine the multiplicity of all its distinct eigenvalues.The following Theorem shows that how we can find all distinct eigenvalues of a class of distance regular graphs.
%%%%%%%%%%%%%%%%%%%%%%%%%%%%%%%%%%%%%%%%%%%%%%%%%%%%%%%%%%%%%%%%%%%%%%%%%%%%%%%%%

\begin{theorem}  \label{b.3}
	{\rm(\cite{paper6}, chapter 4)}
	Let $\Gamma$ be a distance regular graph with  valency  $k$, diameter $d$, adjacency matrix $A$, and intersection array
	$$\{b_0, b_1, ... , b_{d-1}; c_1, c_2, ... , c_d\}.$$
	Then, the tridiagonal $(d+1)\times(d+1)$ matrix
	\begin{center}
		$	B=\begin{bmatrix}
		a_0 & b_0 &0 &0 &...& \\ 
		c_1 & a_1 & b_1&0&...& \\
		0& c_2 & a_2&b_2&0 &... \\
		& & &...& \\
		& & & & \\
		0&... & &0 &c_{d-2}&a_{d-2}&b_{d-2}&0 \\
		0& ...& &0&0&c_{d-1}&a_{d-1}&b_{d-1} \\
		0& ...& &0&0&0&c_d&a_d \\
		\end{bmatrix},$
	\end{center}
	determines all the eigenvalues of $\Gamma$.
\end{theorem}
%%%%%%%%%%%%%%%%%%%%%%%%%%%%%%%%%%%%%%%%%%%%%%%%%%%%%%%%%%%%%%%%%%%%%%%%%%%%%%%%%

\begin{proposition}\label{c.5}
	The bipartite double graph $2.O_k$ of the Odd graph $O_k$ is a cover for the Odd graph $O_k$.
\end{proposition}
\begin{proof}
	The bipartite double graph  $2.O_k$ of the Odd graph $O_k$ has the property that, for each vertex $(v, i)$, there is a
	unique vertex in $2.O_k$  at distance $2k-1$,  from $(v, i)$, say $(v, i^c)$. Thus $V(2.O_k)$ can partitioned into
	$\frac{1}{2}\binom{2k}{k}$ pairs, and these pairs are the fibres of a covering map from $2.O_k$
	onto $O_k$. \qed
\end{proof}
%%%%%%%%%%%%%%%%%%%%%%%%%%%%%%%%%%%%%%%%%%%%%%%%%%%%%%%%%%%%%%%%%%%%%%%%%%%%%%%%%
\begin{proposition}\label{c.6}
	All the eigenvalues of the bipartite double graph $2.O_k$ of the Odd graph $O_k$ are the integers $\lambda_i=\pm(k-i)$ with multiplicity $m(\lambda_i)=\binom{2k-1}{i}-\binom{2k-1}{i-1}$, for $0\leq  i\leq k-1$.
\end{proposition}
	\begin{proof}
		We know that all the eigenvalues of the Odd graph $O_k$  are the integers $\lambda_i = (- 1)^i(k - i)$ with multiplicity
		$ m(\lambda_i)=\binom{2k-1}{i}-\binom{2k-1}{i-1}$ for $0\leq  i\leq k-1$  (see \cite{paper5} Page 74), and it is well known that the eigenvalues of the base graph are also eigenvalues of the cover. Therefore,  all the eigenvalues of $2.O_k$ are the integers $\lambda_i=\pm(k-i)$ with multiplicity $m(\lambda_i)=\binom{2k-1}{i}-\binom{2k-1}{i-1}$ for $0\leq  i\leq k-1$, because the bipartite double graph $2.O_k$ of the Odd graph $O_k$ is a bipartite cover graph over $O_k$. \qed
	\end{proof}
%%%%%%%%%%%%%%%%%%%%%%%%%%%%%%%%%%%%%%%%%%%%%%%%%%%%%%%%%%%%%%%%%%%%%%%%%%%%%%%%%

\begin{lemma}\label{c.7}
	Let $k$ be a positive integer with $k\geq 2$. Then, all the eigenvalues of
	the tridiagonal $(2k)\times(2k)$ matrix
	
	\begin{center}
		$B=\begin{bmatrix}
		0 & k &0 &0 &...& \\ 
		1 & 0 & k-1&0&...& \\
		0& 1 & 0&k-1&0 &... \\
		& & &...& \\
		& & & & \\
		0&... & &0 &k-1&0&1&0 \\
		0& ...& &0&0&k-1&0&1 \\
		0& ...& &0&0&0&k&0 \\
		\end{bmatrix},$
	\end{center}
	are the integers $\lambda_i=\pm(k-i)$  for $0\leq  i\leq k-1$.
	\end{lemma}
\begin{proof}
		Based on (\cite{paper3} Page 201),  the bipartite double graph $2.O_k$ of the Odd graph $O_k$ is a
		distance regular graph, such that the valency of each vertex is $k$, with diameter $2k-1$, and whose intersection array is
		\begin{center}
			$\iota(\Gamma)=\begin{Bmatrix}\ast&1&1&2&2&...&k-1&k-1&k\\
			0&0&0&0&0&...&0&0&0\\k&k-1&k-1&k-2&k-2&...&1&1&\ast\end{Bmatrix}.$
		\end{center}
		Then, by Theorem \ref{b.3}, the tridiagonal $(2k)\times(2k)$ matrix $B$
		determines all the eigenvalues of $2.O_k$. On the other hand, by Proposition \ref{c.6}, all the eigenvalues of $2.O_k$ are the integers $\lambda_i=\pm(k-i)$. Therefore, all the eigenvalues of the tridiagonal $(2k)\times(2k)$  matrix $B$ are the integers $\lambda_i=\pm(k-i)$ for $0\leq  i\leq k-1$. \qed
	\end{proof}
%%%%%%%%%%%%%%%%%%%%%%%%%%%%%%%%%%%%%%%%%%%%%%%%%%%%%%%%%%%%%%%%%%%%%%%%%%%%%%%%%

The following result also appeared in \cite{paper6b}, but our approach for proving it is different.
\begin{theorem}\label{c.8}
\begin{enumerate} 
	\item[(i)]
	 The extended graph $E(2.O_k)$ is an integral.
		
	\item[(ii)]
	  If $k$ is an odd integer then all the eigenvalues of $E(2.O_k)$ are 
	$\theta^{+}_i= +(k-i+(-1)^i)$, and $\theta^{-}_i= -(k-i+(-1)^i)$ for $0\leq i\leq k-1$ and  the multiplicity $m(\theta^{\pm}_i)$ is $\binom{2k-1}{i}-\binom{2k-1}{i-1}$, for $0\leq  i\leq k-1$.
	
	\item[(iii)]
	 If $k$ is an even integer then all the eigenvalues of $E(2.O_k)$ are 
	$\theta^{+}_i= +(k-i+(-1)^i)$, and $\theta^{-}_i= -(k-i+(-1)^i)$ for $0\leq i\leq k-2$ and  the multiplicity $m(\theta^{\pm}_i)$ is $\binom{2k-1}{i}-\binom{2k-1}{i-1}$, for $0\leq  i\leq k-2$, and $m(\theta^{\pm}_{k-1})=2(\binom{2k-1}{k-1}-\binom{2k-1}{k-2})$.  
	\end{enumerate} 
	\end{theorem}
	\begin{proof}
\begin{enumerate} 
		\item[(i)]
		 It is well known that the extended graph $E(2.O_k)$ is an integral graph, see \cite{paper6b}. Our main objective is to determine all eigenvalues of the extended graph $E(2.O_k)$ by leveraging the relationship between its equitable and orbit partitions. First, we show that  the extended graph $E(2.O_k)$ is an integral graph. For this purpose,
		let $\Lambda=E(2.O_k)$ be the extended graph of the bipartite double graph $2.O_k$, is defined already. Based on Theorem \ref{a.3}, and Proposition \ref{c.3},
		for any vertex $(v, i)$ of $E(2.O_k)$ the  distance partition 
		$\Lambda_r((v, i))=\{\Lambda_0(v, i), ... , \Lambda_{k}(v, i)\}$, where $r$ is a non-negative integer not exceeding $k$, the diameter of $E(2.O_k)$ cannot be an equitable partition of $E(2.O_k)$, because  $E(2.O_k)$ cannot be a distance regular graph. 
		Now, let $\Gamma=2.O_k$ be the bipartite double of the Odd graph $O_k$.
		Based on Theorem \ref{a.3}, for any vertex $(v, i)$ of $2.O_k$ the distance partition 
		$\Gamma_r (v,i)=\{\Gamma_0(v, i), ... , \Gamma_{2k-1}(v, i)\}$,  where $r$ is a non-negative integer not exceeding $2k-1$, the diameter of $2.O_k$
		is an equitable partition for $2.O_k$, because $2.O_k$ is a distance regular graph.
		On the other hand, $Aut(2.O_k)$ acts distance transitively on the vertex set of  $2.O_k$, because $2.O_k$ is a distance transitive graph, that is for any vertex $(v, i)$ of $2.O_k$ if $d_\Gamma((v,i), (u,i))=d_\Gamma((v,i), (w,i))=r$, then there is an automorphism $\theta$ of $(Aut(\Gamma))_{(v,i)}$ so that $\theta (u, i)=(w,i)$, and hence $Aut(\Gamma)$ acts transitively on $\Gamma_r (v,i)$, for $0\leq  r\leq 2k-1$.  Therefore, the cells of the distance partition 
		$\Gamma_r (v,i)=\{\Gamma_0(v, i), ... , \Gamma_{2k-1}(v, i)\}$ are the orbits of $(Aut(\Gamma))_{(v,i)}$. 
		Especially, $\Gamma_r (v,i)=\{\Gamma_0(v, i), ... , \Gamma_{2k-1}(v, i)\}$ is an orbit partition of $E(2.O_k)$, because $E(2.O_k)$ is also a vertex transitive graph so that $Aut(E(2.O_k))=Aut(2.O_k)$.
		It is easy  to show that the partition matrix $(B_{\Gamma_r (v,i)})_{\frac{\Gamma}{\Gamma_r (v,i)}}$ of $\Gamma$ is equal to the tridiagonal $(2k)\times(2k)$ matrix
		
		\begin{center}
			$(B_{\Gamma_r (v,i)})_{\frac{\Gamma}{\Gamma_r (v,i)}}=\begin{bmatrix}
			0 & k &0 &0 &...& \\ 
			1 & 0 & k-1&0&...& \\
			0& 1 & 0&k-1&0 &... \\
			& & &...& \\
			& & & & \\
			0&... & &0 &k-1&0&1&0 \\
			0& ...& &0&0&k-1&0&1 \\
			0& ...& &0&0&0&k&0 \\
			\end{bmatrix}.$
		\end{center}
		So, by  Lemma \ref{c.7},  all the eigenvalues of $(B_{\Gamma_r (v,i)})_{\frac{\Gamma}{\Gamma_r (v,i)}}$ are the integers $\lambda_i=\pm (k-i)$.
		On the other hand, the partition matrix $(B_{\Gamma_r (v,i)})_{\frac{\Lambda}{\Gamma_r (v,i)}}$ of graph $\Lambda$ is equal
		to  the  $(2k)\times(2k)$ matrix,
		
		\begin{center}
			$(B_{\Gamma_r (v,i)})_{\frac{\Lambda}{\Gamma_r (v,i)}}=\begin{bmatrix}
			0 & k &0 &0 &...&0&0&0&1 \\
			1 & 0 & k-1&0&...&0&0&1&0 \\
			0& 1 & 0&k-1&0&...&1&0&0 \\
			&  & &...& \\
			&  & &  & \\
			0& 0 &1&0 &...&k-1 &0&1&0 \\
			0&1&0&0  & ...&0&k-1&0&1 \\
			1&0&0& 0 & ...& 0&0&k&0 \\
			\end{bmatrix}.$
		\end{center}
		Hence, we have
		$(B_{\Gamma_r (v,i)})_{\frac{\Lambda}{\Gamma_r (v,i)}}=(B_{\Gamma_r (v,i)})_{\frac{\Gamma}{\Gamma_r (v,i)}}+C$, where $C$ is equal
		to the  $(2k)\times(2k)$ matrix,
		\begin{center}
			$C=\begin{bmatrix}
			0 & 0 &0 &0 &...&0&0&0&1 \\
			0 & 0 & 0&0&...&0&0&1&0 \\
			0& 0 & 0&0& &0&1&0&0 \\
			&  & &...& \\
			&  & &  & \\
			0& 0 &1&0 &...&0 &0&0&0 \\
			0&1&0&0  & ...&0&0&0&0 \\
			1&0&0& 0 & ...& 0&0&0&0 \\
			\end{bmatrix}.$
		\end{center}
		Since, the  $(2k)\times(2k)$ matrix  $(B_{\Gamma_r (v,i)})_{\frac{\Gamma}{\Gamma_r (v,i)}}$ has $2k$ distinct eigenvalues, then $(B_{\Gamma_r (v,i)})_{\frac{\Gamma}{\Gamma_r (v,i)}}$ is diagonalizable. Aalso, we can verify that the $(2k)\times(2k)$ matrix $C$ is diagonalizable, and
		$(B_{\Gamma_r (v,i)})_{\frac{\Gamma}{\Gamma_r (v,i)}}C=C(B_{\Gamma_r (v,i)})_{\frac{\Gamma}{\Gamma_r (v,i)}}$. 
		Therefore, $(B_{\Gamma_r (v,i)})_{\frac{\Gamma}{\Gamma_r (v,i)}}$ and $C$ are  simultaneously diagonalizable, and hence there is a  base $W=\{w_1, w_2, ... , w_{2k}\}$ of $\mathbb{R}^{2k}$, such that each $w_i\in W$ is an eigenvector of $(B_{\Gamma_r (v,i)})_{\frac{\Gamma}{\Gamma_r (v,i)}}$ and $C$.
		Since, 
		$$(B_{\Gamma_r (v,i)})_{\frac{\Gamma}{\Gamma_r (v,i)}}w_i+Cw_i=\lambda_iw_i\pm 1w_i=(\lambda_i\pm 1)w_i,$$
		then each $w_i\in W$ is an eigenvector of  $(B_{\Gamma_r (v,i)})_{\frac{\Lambda}{\Gamma_r (v,i)}}$, because $(B_{\Gamma_r (v,i)})_{\frac{\Lambda}{\Gamma_r (v,i)}}=(B_{\Gamma_r (v,i)})_{\frac{\Gamma}{\Gamma_r (v,i)}}+C$. It follows that, all the eigenvalues of $(B_{\Gamma_r (v,i)})_{\frac{\Lambda}{\Gamma_r (v,i)}}$ are integers.
		Hence $E(2.O_k)$ is an integral graph, because based on Theorem \ref{b.6},
		each eigenvalue of $E(2.O_k)$ is an eigenvalue of $(B_{\Gamma_r (v,i)})_{\frac{\Lambda}{\Gamma_r (v,i)}}$.

		\item[(ii)] 
		Let $k$ be an odd integer, $\Gamma=2.O_k$ and $\Lambda=E(2.O_k)$.  Now, let $\lambda^{+}_i= +(k-i)$ and $\lambda^{-}_i= -(k-i)$ denote all the eigenvalues of the
 Odd graph $O_k$ with the multiplicity $m(\lambda^{\pm}_i)=\binom{2k-1}{i}-\binom{2k-1}{i-1}$, for $0\leq  i\leq k-1$.  First, we show that all the eigenvalues of  the extended graph $E(2.O_k)$ are $\theta^{+}_i= +(k-i+(-1)^i)$, and $\theta^{-}_i= -(k-i+(-1)^i)$ for $0\leq i\leq k-1$,
		and then we show that the multiplicity $m(\theta^{\pm}_i)=\binom{2k-1}{i}-\binom{2k-1}{i-1}$, for $0\leq  i\leq k-1$. 
		For this purpose, let $A_\Gamma$ and $A_\Lambda$ be the adjacency matrices of graphs $2.O_k$ and $E(2.O_k)$, respectively. Hence, 
		$A_\Lambda=A_\Gamma+ D$, where $D$ is equal to the  $(q)\times(q)$ matrix
		\begin{center}
			$D=\begin{bmatrix}
			0 & 0 &0 &0 &...&0&0&0&1 \\
			0 & 0 & 0&0&...&0&0&1&0 \\
			0& 0 & 0&0& &0&1&0&0 \\
			&  & &...& \\
			&  & &  & \\
			0& 0 &1&0 &...&0 &0&0&0 \\
			0&1&0&0  & ...&0&0&0&0 \\
			1&0&0& 0 & ...& 0&0&0&0 \\
			\end{bmatrix},$
		\end{center}
		where $q=2\binom{2k-1}{k-1}$.
		It is well known that, $2.O_k$ is a distance regular graph with the diameter $2k-1$, hence $2.O_k$ has $2k$ eigenspace $V^{(\pm)}_0, ..., V^{(\pm)}_{k-1}$ so that $dim(V^{(\pm)}_i)=\binom{2k-1}{i}-\binom{2k-1}{i-1}$, for $0\leq  i\leq k-1$. In particular, we can verify that $A_\Gamma$ and $D$ are  commute and  diagonalizable, hence they are simultaneously diagonalizable, and hence each eigenspace $V^{(\pm)}_i$ is invariant on $D$. Now, we show that $D$ acts as a scalar on each $V^{(\pm)}_i$.
For this purpose, let $A_{O_k}$  be the adjacency matrice of the Odd graph $O_k$ and $E_i$ be the $i^{th}$ eigenspace of $O_k$ (for eigenvalue $\mu_i=(-1)^i (k-i)$,  $0\leq  i\leq k-1$). Hence,  if $\mu_i$ is a positive eigenvalue, then  eigenspace $V_i^{+}$ corresponds to eigenvectors of 
form $(\phi, \phi)$, where $\phi \in E_i$, and similarly if $\mu_i$ is a negative eigenvalue, then  eigenspace $V_i^{+}$ corresponds to eigenvectors of 
form $(\phi, -\phi)$. Therefore, if $A_{O_k}\phi=\mu_i\phi$ and $\mu_i$ positive eigenvalue,  then $V_i^{+}$ has eigenvectors $(\phi, \phi)$ with $\phi \in E_i$, and $V_i^{-}$ has $(\phi, -\phi)$ with $\phi \in E_i$,
and similarly if $\mu_i$ is a negative eigenvalue, then $V_i^{+}$ has eigenvectors $(\phi, -\phi)$ with $\phi \in E_i$, and $V_i^{-}$ has $(\phi, \phi)$ with $\phi \in E_i$. Thus,
$D$ acts as a scalar on each $V^{(\pm)}_i$, specifically $D|_{V^{(+)}_i} = (-1)^i \cdot I$  and $D|_{V^{(-)}_i} = -(-1)^i \cdot I$ for $0\leq  i\leq k-1$. 
That is, every vector in $V^{(+)}_i$ is an eigenvector of $D$ with eigenvalue $(-1)^i$ and every vector in $V^{(-)}_i$ is an eigenvector of $D$ with eigenvalue $-(-1)^i$.
Therefore the eigenvalues of $D$ on the $i^{th}$ eigenspace of  $A_\Gamma$ are $(-1)^i$ and  $-(-1)^i$.
Now, let $Z=\{z_1, z_2, ... , z_{q}\}$  such that each $z_i\in Z$ be an eigenvector of $A_\Gamma$ and $D$ . Hence, if $z_i\in V^{(+)}_i$ then 
$A_\Lambda z_i=A_\Gamma z_i+ D z_i=+(k-i+(-1)^i)$, and similarly, for  $z_i\in V^{(-)}_i$ we have $A_\Lambda z_i=-(k-i+(-1)^i)$.
 Since $D$ acts as a scalar on this eigenspace, then the multiplicity $m(\theta^{\pm}_i)=\binom{2k-1}{i}-\binom{2k-1}{i-1}$, for $0\leq  i\leq k-1$. 
		
				\item[(iii)]
		 Now, let $k$ be an even integer. By similar way is done in the proof of (ii), we can verify that if $k$ is an even integer
		then  the eigenvalues of $E(2.O_k)$ are 
		$\theta^{+}_i= +(k-i+(-1)^i)$, and $\theta^{-}_i= -(k-i+(-1)^i)$ for $0\leq i\leq k-2$, each with  multiplicity is $\binom{2k-1}{i}-\binom{2k-1}{i-1}$, for $0\leq  i\leq k-2$. Note that, if $k$ is even integer and  $i=k-1$, then $\theta^{\pm}_i= \pm 0$, 
		and hence $m(\theta^{\pm}_{k-1})=2(\binom{2k-1}{k-1}-\binom{2k-1}{k-2})$.  \qed
		\end{enumerate} 
	\end{proof}
%%%%%%%%%%%%%%%%%%%%%%%%%%%%%%%%%%%%%%%%%%%%%%%%%%%%%%%%%%%%%%%%%%%%%%%%%%%%%%%%%

\begin{theorem}\label{c.9}{\rm\cite{paper6-dd}}
	The spectrum of Johnson graph $J(n,m)$, $n,m\in\mathbb{N}$, $m<n$, is
	$$
	\{\lambda_{i}^{m_{i}}|\ \lambda_{i}=(m-i)(n-m-i)-i,\ m({\lambda_i})=\binom{n}{i}-\binom{n}{i-1},\ 0\leq i\leq d\},
	$$
	where \(d=\min(m,n-m)\) is the diameter of $J(n,m)$. 
\end{theorem}
%%%%%%%%%%%%%%%%%%%%%%%%%%%%%%%%%%%%%%%%%%%%%%%%%%%%%%%%%%%%%%%%%%%%%%%%%%%%%%%%%
 Although the eigenvalues of the enhanced Johnson graph $EJ(2m,m)$ is computed  in Theorem 3.9 of \cite{w-3}, but the multiplicity of these distinct eigenvalues has not yet been calculated. Hence, in the following result we obtain all the eigenvalues of the enhanced Johnson graph $EJ(2m,m)$ by another approach and we use our techniques to fully determine the multiplicity of all distinct eigenvalues for $EJ(2m,m)$.

\begin{theorem}\label{c.11}
	If $m$ is an odd positive integer is greater than or equal to $3$, and $\theta_i$ is an eigenvalue of  the enhanced Johnson graph $EJ(2m,m)$, then the spectrum of enhanced Johnson graph $EJ(2m,m)$ is
	$$
	\{\theta_{i}^{m_{i}}|\ \theta_{i}=\lambda_{i}+(-1)^i,\ m({\theta_i})=\binom{2m}{i}-\binom{2m}{i-1},\ 0\leq i\leq m\},
	$$
	where $\lambda_{i}$ is an eigenvalue of  the Johnson graph $J(2m,m)$.
	\end{theorem}
	\begin{proof}
		Let $\Gamma=J(2m,m)$ and $\Lambda=EJ(2m,m)$.  First we show that all the eigenvalues of  the enhanced Johnson graph $EJ(2m,m)$ are $\theta_{i}=\lambda_{i}+(-1)^i$, for $0\leq i\leq m$,
		and then we show that the multiplicity $m(\theta_i)=\binom{2m}{i}-\binom{2m}{i-1}$, for $0\leq  i\leq m$. 
		For this purpose, let $A_\Gamma$ and $A_\Lambda$ be the adjacency matrices of graphs $J(2m,m)$ and $EJ(2m,m)$, respectively. Hence, 
		$A_\Lambda=A_\Gamma+ D$, where $D$ is equal to the  $(q)\times(q)$ matrix
		\begin{center}
			$D=\begin{bmatrix}
			0 & 0 &0 &0 &...&0&0&0&1 \\
			0 & 0 & 0&0&...&0&0&1&0 \\
			0& 0 & 0&0& &0&1&0&0 \\
			&  & &...& \\
			&  & &  & \\
			0& 0 &1&0 &...&0 &0&0&0 \\
			0&1&0&0  & ...&0&0&0&0 \\
			1&0&0& 0 & ...& 0&0&0&0 \\
			\end{bmatrix},$
		\end{center}
		where $q=\binom{2m}{m}$.
		It is well known that, $J(2m,m)$ is a distance regular graph with the diameter $m$, hence $J(2m,m)$ has $m+1$ eigenspace $V_0, ..., V_m$ so that $dim(V_i)=\binom{2m}{i}-\binom{2m}{i-1}$, for $0\leq  i\leq m$. Based on \cite{bbb2}, $J(2m,m)$ belongs to a symmetric association scheme generated by distance adjacency matrices $A_0, A_1, \dots, A_m$, then $A_m=D$. Therefore, $A_\Gamma$ and $D$ are  commute and  diagonalizable, hence they are simultaneously diagonalizable, and hence each eigenspace $V_i$ is invariant on $D$.
Also, based on (\cite{paper67a} Page 48), the eigenvalues of $A_m=D$ on $V_i$ in $J(2m,m)$ are given by Eberlein polynomial

$$
E_m(i) = \sum_{j=0}^{m} (-1)^j \binom{i}{j} \binom{m-i}{m-j}^2.
$$
We can verify that, $\binom{i}{j}=0$ if $j>i$ and  $\binom{m-i}{m-j}=0$ if $i>j$. Hence, if $i=j$, then we have $E_m(i)=(-1)^i$. So $D$ acts as a scalar on each $V_i$, specifically $D|_{V_i} = (-1)^i \cdot I$.
Thus the eigenvalues of $D$ on the $i^{th}$ eigenspace of  $A_\Gamma$ are $(-1)^i$.
Now, let $Z=\{z_1, z_2, ... , z_{q}\}$  such that each $z_i\in Z$ be an eigenvector of $A_\Gamma$ and $D$ . Hence, for any eigenvector $z_i\in Z$ in the $i^{th}$ eigenspace of $A_\Gamma$ we have $A_\Gamma z_i=\lambda_iz_i$ and $Dz_i=(-1)^iz_i$, therefore 
		$A_\Lambda z_i=(\lambda_i+(-1)^i)z_i$. Since $D$ acts as a scalar on this eigenspace, then the multiplicity $m(\theta_i)=\binom{2m}{i}-\binom{2m}{i-1}$, for $0\leq  i\leq m$. \qed
	\end{proof}
%%%%%%%%%%%%%%%%%%%%%%%%%%%%%%%%%%%%%%%%%%%%%%%%%%%%%%%%%%%%%%%%%%%%%%%%%%%%%%%%%

\begin{theorem}\label{c.12}
	If $m$ is an even positive integer is greater than or equal to $2$, and $\theta_i$ is an eigenvalue of  the enhanced Johnson graph $EJ(2m,m)$, then the spectrum of enhanced 
	Johnson graph $EJ(2m,m)$ is
	$$
	\{\theta_{i}^{m_{i}}|\ \theta_{i}=\lambda_{i}+(-1)^i,\ m({\theta_i})=\binom{2m}{i}-\binom{2m}{i-1},\ 0\leq i\leq m-2\},
	$$
	also $m({\theta_{m-1}})=m({\lambda_{m-1}})+m({\lambda_{m}})$, where $\lambda_{i}$ is an eigenvalue of  the Johnson graph $J(2m,m)$.
	\end{theorem} 
	\begin{proof}
		Let $\Gamma=J(2m,m)$ and $\Lambda=EJ(2m,m)$. Based on the proof of previous Theorem, we can show that all the eigenvalues of  the enhanced Johnson graph $EJ(2m,m)$ is $\theta_{i}=\lambda_{i}+(-1)^i$, for $0\leq i\leq m$,
		and then we can show that the multiplicity $m(\theta_i)=\binom{2m}{i}-\binom{2m}{i-1}$, for $0\leq  i\leq m-2$.
		Note that, if $m$ is an even integer then by according to $\lambda_{i}=(m-i)(2m-m-i)-i$ for $0\leq i\leq m$, we have 
		$\lambda_{m-1} = 2-m$ and  $\lambda_{m} = -m$. Therefore, by according to
		$\theta_{i}=\lambda_{i}+(-1)^i$, for $0\leq i\leq m$, we have 
		$\theta_{m-1}=\theta_{m}=1-m$,
		and hence  $m({\theta_{m-1}})=m({\lambda_{m-1}})+m({\lambda_{m}})$. \qed
		
	\end{proof}
%%%%%%%%%%%%%%%%%%%%%%%%%%%%%%%%%%%%%%%%%%%%%%%%%%%%%%%%%%%%%%%%%%%%%%%%%%%%%%%%%
%%%%%%%%%%%%%%%%%%%%%%%%%%%%%%%%%%%%%%%%%%%%%%%%%%%%%%%%%%%%%%%%%%%%%%%%%%%%%%%%%

\bigskip

\begin{thebibliography}{5}
%%%%%%%%%%%%%%%%%%%%%%%%%%%%%%%%%%%%%%%%%%%%%%%%%%%%%%%%%%%%%%%%%%%%%%%%%%%%%%%%%
\bibitem{paper0}{A. Abdollahi and E. Vatandoost}, Which Cayley graphs are integral?, \emph{Electron. J. Combin.}, \textbf{16(1)} (2009), R122, 1-17.
%%%%%%%%%%%%%%%%%%%%%%%%%%%%%%%%%%%%%%%%%%%%%%%%%%%%%%%%%%%%%%%%%%%%%%%%%%%%%%%%%
\bibitem{paper1}{O. Ahmadi, N. Alon, L. F. Blake and I. E. Shparlinski}, Graphs with integral spectrum, \emph{Linear Alg. Appl.},
\textbf{430} (2009), 547-552.
%%%%%%%%%%%%%%%%%%%%%%%%%%%%%%%%%%%%%%%%%%%%%%%%%%%%%%%%%%%%%%%%%%%%%%%%%%%%%%%%%
\bibitem{paper1.2}{A. T. Balaban, D. Farcussiu and R. Banica}, \emph{Graphs of multiple 1; 2-shifts in carbonium ions
		and related systems}, Rev. Roum. Chim., \textbf{11} (1966), 1205-1227.
%%%%%%%%%%%%%%%%%%%%%%%%%%%%%%%%%%%%%%%%%%%%%%%%%%%%%%%%%%%%%%%%%%%%%%%%%%%%%%%%%
	\bibitem{paper2}{K. Balinska et al.}, A survey on integral graphs,\emph{ Univ. Beograd, Publ. Elektrotehn. Fak. Ser. Mat.},
	\textbf{13} (2003),  42-65.
%%%%%%%%%%%%%%%%%%%%%%%%%%%%%%%%%%%%%%%%%%%%%%%%%%%%%%%%%%%%%%%%%%%%%%%%%%%%%%%%%
\bibitem{bbb2}{E. Bannai, T. Ito}, Algebraic combinatorics. I. The Benjamin/Cummings Publishing Co., Inc., Menlo Park, CA, (1984).
%%%%%%%%%%%%%%%%%%%%%%%%%%%%%%%%%%%%%%%%%%%%%%%%%%%%%%%%%%%%%%%%%%%%%%%%%%%%%%%%%
	\bibitem{paper3}{E. Bannai}, Codes in bipartite distance regular graphs,
	\emph{journal of the london mathematical society},
	\textbf{s2-16 (2)} (1977), 197-202.
%%%%%%%%%%%%%%%%%%%%%%%%%%%%%%%%%%%%%%%%%%%%%%%%%%%%%%%%%%%%%%%%%%%%%%%%%%%%%%%%%
	\bibitem{paper3.1} {N. Biggs}, \emph{An edge coloring problem}, Amer. Math. Montly, \textbf{79} (1972), 1018-1020.
%%%%%%%%%%%%%%%%%%%%%%%%%%%%%%%%%%%%%%%%%%%%%%%%%%%%%%%%%%%%%%%%%%%%%%%%%%%%%%%%%
		\bibitem{paper5}{N. Biggs}, \emph{Some Odd Graph Theory}, \emph{Ann. New York Acad. Sci.},
	\textbf{319} (1979), 71-81.
	%%%%%%%%%%%%%%%%%%%%%%%%%%%%%%%%%%%%%%%%%%%%%%%%%%%%%%%%%%%%%%%%%%%%%%%%%%%%%%%%%
	\bibitem{paper4} {N. Biggs}, Algebraic Graph Theory, Cambridge University Press, Cambridge (1974), second
	edition (1993).
%%%%%%%%%%%%%%%%%%%%%%%%%%%%%%%%%%%%%%%%%%%%%%%%%%%%%%%%%%%%%%%%%%%%%%%%%%%%%%%%%
	\bibitem{paper6}{A. E. Brower, A. M. Cohen and  A. Neumaier}, \emph{Distance Regular Graphs}, Springer, Berlin, (1989).
	%%%%%%%%%%%%%%%%%%%%%%%%%%%%%%%%%%%%%%%%%%%%%%%%%%%%%%%%%%%%%%%%%%%%%%%%%%%%%%%%%
	\bibitem{paper6-dd}{A. E. Brower and W. H Haemers},  \emph{Spectra of Graphs}, Springer, (2012).
	%%%%%%%%%%%%%%%%%%%%%%%%%%%%%%%%%%%%%%%%%%%%%%%%%%%%%%%%%%%%%%%%%%%%%%%%%%%%%%%%%
	\bibitem{paper6a}{C. Dalf\'{o}, M. A. Fiol and M. Mitjana}, \emph{On Middle Cube Graphs}, Electronic Journal of Graph Theory and Applications, \textbf{3 (2)} (2015), 133-145.
	%%%%%%%%%%%%%%%%%%%%%%%%%%%%%%%%%%%%%%%%%%%%%%%%%%%%%%%%%%%%%%%%%%%%%%%%%%%%%%%%%
\bibitem{paper67a}{P. Delsarte},  An algebraic approach to the association schemes of coding theory, Centrex Publishing Co., (1973). 
	%%%%%%%%%%%%%%%%%%%%%%%%%%%%%%%%%%%%%%%%%%%%%%%%%%%%%%%%%%%%%%%%%%%%%%%%%%%%%%%%%
	\bibitem{paper6b}{M. A. Fiol, E. Garriga and J. L. A. Yebra}, \emph{On twisted odd graphs}, Combin. Probab. Comput.,
	\textbf{9} (2000), 227-240.
	%%%%%%%%%%%%%%%%%%%%%%%%%%%%%%%%%%%%%%%%%%%%%%%%%%%%%%%%%%%%%%%%%%%%%%%%%%%%%%%%%
	\bibitem{paper7}{C. Godsil and   G. Royle}, \emph{Algebraic Graph Theory}, Springer, New York, (2001).
%%%%%%%%%%%%%%%%%%%%%%%%%%%%%%%%%%%%%%%%%%%%%%%%%%%%%%%%%%%%%%%%%%%%%%%%%%%%%%%%%
	\bibitem{paper8}{F. Harary and A. Schwenk}, Which graphs have integral spectra?, Lect.
	Notes Math., \emph{Springer Verlag},
	\textbf{406} (1974), 45-50.
%%%%%%%%%%%%%%%%%%%%%%%%%%%%%%%%%%%%%%%%%%%%%%%%%%%%%%%%%%%%%%%%%%%%%%%%%%%%%%%%%
\bibitem{paper1Ee-1}{J. Kim, S. Kim, E. Cheng and László Lipták}, Topological properties of folded hyper-star networks, \emph{J Supercomput}, \textbf{59} (2012), 1336–1347.
	%%%%%%%%%%%%%%%%%%%%%%%%%%%%%%%%%%%%%%%%%%%%%%%%%%%%%%%%%%%%%%%%%%%%%%%%%%%%%%%%%
			\bibitem{paper9}{S. Kudose}, Equitable partition and orbit partition,
	\emph{The University of Chicago Department of Mathematics}, 1-9.
	%%%%%%%%%%%%%%%%%%%%%%%%%%%%%%%%%%%%%%%%%%%%%%%%%%%%%%%%%%%%%%%%%%%%%%%%%%%%%%%%%
		\bibitem{w-1}{S. M. Mirafzal and  A. Zafari}, \emph{Some algebraic properties of bipartite Kneser
		graphs}, Ars Combin., \textbf{153} (2020), 3–12.
%%%%%%%%%%%%%%%%%%%%%%%%%%%%%%%%%%%%%%%%%%%%%%%%%%%%%%%%%%%%%%%%%%%%%%%%%%%%%%%%%
	\bibitem{w-2}{S. M. Mirafzal}, \emph{A new class of integral graphs constructed from the hypercube}, Linear Algebra
	Appl., \textbf{558}  (2018), 186-194.
	%%%%%%%%%%%%%%%%%%%%%%%%%%%%%%%%%%%%%%%%%%%%%%%%%%%%%%%%%%%%%%%%%%%%%%%%%%%%%%%%%
	\bibitem{w-3}{S. M. Mirafzal and M. Ziaee}, \emph{Some algebraic aspects of enhanced Johnson graphs}, Acta Math. Univ. Comenianae
	Vol. LXXXVIII, \textbf{2} (2019), 257-266.
	%%%%%%%%%%%%%%%%%%%%%%%%%%%%%%%%%%%%%%%%%%%%%%%%%%%%%%%%%%%%%%%%%%%%%%%%%%%%%%%%%
		\bibitem{paper1E-1}{A. Staples-Moore}, Equitable partitions in graph theory,
	\emph{The University of Chicago Department of Mathematics}, 1-10.
	
	\end{thebibliography}
\end{document}